\documentclass[11pt]{article}

  \usepackage[letterpaper, left=1in, right=1in, top=1in, bottom=1in]{geometry}
  
  \usepackage[affil-it]{authblk}  
  
  \usepackage{palatino}  
  \usepackage{setspace}  
  \usepackage{color}  
  \usepackage[dvipsnames]{xcolor}  

	\usepackage[labelfont=bf]{caption}	

  \usepackage{amsmath}  
  \usepackage{amssymb}  
  \usepackage{amsthm}  
  \allowdisplaybreaks[4]  
	\usepackage{relsize}  
	\usepackage{algorithm}  
	\usepackage{algpseudocode}  
  
  \DeclareMathOperator*{\minimize}{minimize}
\DeclareMathOperator*{\maximize}{maximize}
  \DeclareMathOperator*{\argmin}{arg\,min}
  \newcommand{\st}{\mathrm{subject\;to}}

	\theoremstyle{definition}  

        \newtheorem{remark}{Remark}
	
	\usepackage{longtable}  
	\usepackage{multirow}  
	\usepackage{array}  
	\usepackage{booktabs}  
    
\usepackage{graphicx}  
\usepackage{float}  
\usepackage{subfig}  
\usepackage{pgf,tikz,pgfplots}
\pgfplotsset{compat=1.15}
\usepackage{mathrsfs}
\usepackage{epsfig}
\usepackage{epstopdf}
\usepackage{pgfplotstable}
\usepackage{ifthen}

	\usepackage{hyperref}
	\hypersetup{
		pdfstartview = {FitV},  
		colorlinks = true,    
		linkcolor = NavyBlue,  
		citecolor = NavyBlue,  
		filecolor = NavyBlue,  
		urlcolor = NavyBlue  
	}  
 
\makeatletter
\newenvironment{varsubequations}[1]
 {%
  \addtocounter{equation}{-1}%
  \begin{subequations}
  \renewcommand{\theparentequation}{#1}%
  \def\@currentlabel{#1}%
 }
 {%
  \end{subequations}\ignorespacesafterend
 }
 
\makeatother

  \usepackage[numbers, sort&compress, super]{natbib}  

\title{Data-driven decision-focused surrogate modeling}
\author[1]{Rishabh Gupta}
\author[1]{Qi Zhang \thanks{Corresponding author (qizh@umn.edu)}}
\affil[1]{Department of Chemical Engineering and Materials Science, \break University of Minnesota, Minneapolis, MN 55455, USA}
\date{}

\begin{document}

\maketitle

\begin{abstract}
We introduce the concept of decision-focused surrogate modeling for solving computationally challenging nonlinear optimization problems in real-time settings. The proposed data-driven framework seeks to learn a simpler, e.g. convex, surrogate optimization model that is trained to minimize the \emph{decision prediction error}, which is defined as the difference between the optimal solutions of the original and the surrogate optimization models. The learning problem, formulated as a bilevel program, can be viewed as a data-driven inverse optimization problem to which we apply a decomposition-based solution algorithm from previous work. We validate our framework through numerical experiments involving the optimization of common nonlinear chemical processes such as chemical reactors, heat exchanger networks, and material blending systems. We also present a detailed comparison of decision-focused surrogate modeling with standard data-driven surrogate modeling methods and demonstrate that our approach is significantly more data-efficient while producing simple surrogate models with high decision prediction accuracy.
\end{abstract}


\section{Introduction}
\label{sec:Intro}

Efficient and safe process operations require decision making in real time; this is more important than ever as the chemical industry faces new fast-changing markets, greater feedstock variability, and increasingly time-sensitive availability of resources such as intermittent renewable energy. Many online decision-making frameworks, including model predictive control (MPC) and real-time optimization (RTO), involve the solving of mathematical optimization problems. Often, the computational complexity of the optimization problem presents a major challenge such that the long solution time renders it ineffective in online applications. A common approach to tackling this challenge is to perform the online optimization with a \textit{surrogate model}, which is an approximation of the original model that can be solved more efficiently.

Typically in surrogate modeling, one tries to replace complicating functions that are embedded in the optimization problem with simpler ones. This approach is widely used in process systems engineering (PSE) where the complicating functions are often associated with process units that exhibit complex nonlinear behavior. Over the years, the PSE community has developed a myriad of surrogate process models that include shortcut and lumped models derived from first principles and engineering assumptions, reduced-order dynamic models constructed using model reduction methods such as singular perturbation analysis, data-driven models based on, for example, Gaussian processes and artificial neural networks, and many more. Here, a major underlying assumption is that a surrogate model that provides a good approximation of the original functions will, once incorporated in the optimization problem, also lead to solutions that are close to the true optimal solutions. However, it is unclear whether or under what conditions this assumption holds and how accurate the surrogate model needs to be. In fact, as we show in this work, one can easily find examples in which the original optimization problem and the surrogate optimization problem achieve very different optimal solutions despite having a highly accurate embedded surrogate model.

In this work, we propose a new surrogate modeling framework that explicitly aims to construct surrogate models that minimize the \textit{decision prediction error} defined as the difference between the optimal solutions of the original and the surrogate optimization problems; we hence refer to it as \textit{decision-focused surrogate modeling}. We take a data-driven approach in which the original optimization problem is solved offline with different model inputs, resulting in a dataset where each data point is an input-decision pair. We then develop an inverse optimization \cite{Ahuja2001} approach that directly learns from the given data a surrogate optimization model that has a simpler (e.g. convex) form and minimizes the decision prediction error subject to the restrictions on the form of the surrogate model. Results from multiple computational case studies show that the proposed approach outperforms alternative surrogate modeling approaches in decision accuracy, data efficiency, and/or computational efficiency of the resulting surrogate optimization model.

In the remainder of this paper, we first provide a systematic review of the main existing surrogate modeling frameworks and other related works in Section \ref{sec:Background}. We introduce the concept of decision-focused surrogate modeling, provide the corresponding mathematical problem formulation, and present a solution algorithm in Section \ref{sec:DFSM}. Several numerical case studies based on typical nonlinear chemical engineering systems are presented in Section \ref{sec:case_studies}. Finally, we conclude in Section \ref{sec:conclude}.

\section{Background and related work}
\label{sec:Background}

In this section, we present an overview of the major surrogate modeling frameworks and other related work. We formally describe the different approaches to make clear the main conceptual differences, which will help us motivate the proposed decision-focused surrogate modeling framework and highlight its unique features in Section \ref{sec:DFSM}. We focus on data-driven approaches but make reference to other related methods wherever appropriate. 

Without loss of generality, we assume the original optimization problem to be of the following form:
\begin{equation}
\label{eqn:OriginalOP}
    \begin{aligned}
        \minimize_{x \in \mathcal{X}} \quad & f(x, u) \\
        \st \quad & g(x, u) \leq 0 \\
    \end{aligned}
\end{equation}
where $x$ are the decision variables and $\mathcal{X} \subseteq \mathbb{R}^n$.  We assume that the input parameters $u$ can be chosen from a given set $\mathcal{U} \subseteq \mathbb{R}^p$. The constraints describing the feasible region of \eqref{eqn:OriginalOP} are defined by the functions $g:\mathbb{R}^n \times \mathbb{R}^p \rightarrow \mathbb{R}^{k}$, and $f:\mathbb{R}^n \times \mathbb{R}^p \rightarrow \mathbb{R}$ is the objective function of the problem.

\subsection{Optimization with embedded surrogate models} \label{sec:embedded}
We first describe the traditional surrogate modeling framework as outlined in Section \ref{sec:Intro}, which we call optimization with \textit{embedded} surrogate models. To explain the main idea of this approach, we rewrite \eqref{eqn:OriginalOP} as follows:
\begin{equation}
\label{eqn:OriginalOP_embsurr}
    \begin{aligned}
        \minimize_{x \in \mathcal{X}} \quad & f(x, u) \\
        \st \quad & h(x, u) \leq 0 \\
        & d(x, u) \leq 0,
    \end{aligned}
\end{equation}
where we divide the constraint functions $g$ into two sets of functions $h$ and $d$. This surrogate modeling approach is typically applied when problem \eqref{eqn:OriginalOP_embsurr} is such that its computational complexity mainly stems from the functions $d$; hence, the goal is to replace them with a simpler set of functions $\hat{d}$. Such surrogate models are typically generated using simplifying assumptions based on physical and engineering insights, model order reduction techniques, or data-driven methods \citep{Biegler2014}. In a data-driven approach, one first generates a set of $N$ data points where for each point $i \in \mathcal{I} = \{1,\dots,N\}$, $\bar{x}_i$ and $\bar{u}_i$ are sampled from $\mathcal{X}$ and $\mathcal{U}$, respectively, and the corresponding function evaluations $\bar{y}_i = d(\bar{x}_i,\bar{u}_i)$ is computed. Assuming that $\hat{d}(\cdot)$ can be chosen from the set of functions $\mathcal{D}$, one then solves the following empirical risk minimization problem:
\begin{equation}
    \minimize_{\hat{d}(\cdot) \in \mathcal{D}} \; \frac{1}{N} \sum_{i \in \mathcal{I}} \ell \left( \bar{y}_i, \hat{d}(\bar{x}_i,\bar{u}_i) \right),
\end{equation}
where $\ell(\cdot,\cdot)$ denotes a loss function, e.g. the Euclidean distance, that is a measure of the difference between the observation and the estimate. 

Once we have obtained $\hat{d}$, the surrogate model for $d$, we can formulate the following surrogate optimization model:
\begin{equation}
    \begin{aligned}
        \minimize_{x \in \mathcal{X}} \quad & f(x,u) \\
        \st \quad & h(x,u) \leq 0 \\
        & \hat{d}(x,u) \leq 0.
    \end{aligned}
\end{equation}
Note that we use the term \textit{surrogate model} for the alternative model trained using input-output data from the original model and \textit{surrogate optimization model} for the resulting optimization problem that incorporates the surrogate model. 

There exist a large number of statistical learning methods that are used to train these embedded surrogate models; popular choices in PSE include response surface methods \cite{Jones2001, Jia2009, Boukouvala2010}, kriging (or Gaussian process regression) \cite{Caballero2008, Boukouvala2011, Heo2012}, and deep learning \cite{Henao2011, Lee2018, Schweidtmann2019}. Data-driven approaches are often combined with first-principles modeling, resulting in gray-box models. For many physical systems, gray-box models have proven to perform better in terms of model accuracy and interpretability compared to purely data-driven models; hence, first-principles elements are often incorporated as long as they are not overly complex \cite{Cozad2015, Asprion2017, Boukouvala2017a}. A key challenge in surrogate modeling is the balance between model accuracy and computational efficiency; recent efforts aim at generating surrogate models composed of the simplest functional forms possible, given a desired model accuracy, to facilitate their use in mathematical optimization. A prominent example of such an approach is the ALAMO framework proposed by Sahinidis and coworkers \cite{Cozad2014, Wilson2017}.

The constraints in an optimization problem define the feasible region for the variables of the problem; thus, generating surrogate constraint functions can also be interpreted as creating a surrogate representation of the feasible region. Approaches that are directly derived from this interpretation of surrogate modeling consider datasets that consist of feasible and infeasible points. With that, constructing a surrogate model essentially becomes a binary classification problem. Existing works have performed this kind of surrogate modeling using classical classification approaches such as support vector machines (SVMs) \citep{Basudhar2012, Ibrahim2021}. The intricate shapes of the feasible regions arising in many process systems applications have further motivated the development of more complex geometry-based approaches. \citet{Ierapetritou2001} proposes to approximate a convex feasible region with the convex hull of a set of points sampled at the boundary. For nonconvex models, \citet{Goyal2003} develop an approach in which the feasible region is approximated by subtracting outer polytopes around the infeasible region from an expanded convex hull obtained from simplicial approximation \citep{Goyal2002}. The nonconvex case has also been addressed using shape reconstruction techniques \citep{Banerjee2005}. \citet{Zhang2016b} propose an algorithm for constructing a union of multiple polytopes that approximates the nonconvex feasible region from which the data points are sampled. In a similar spirit, \citet{Schweidtmann2022} apply persistent homology, a topological data analysis method, to identify potential holes and clusters in the data; this information is then used to obtain a representation of the nonconvex feasible region using one-class SVMs. In another approach, feasibility is captured using a so-called feasibility function \citep{Halemane1983}, and data-driven approaches are applied to approximate that function; see \citet{Bhosekar2018} for a comprehensive review of this type of methods.

\subsection{Learning optimization proxies}
\label{sec:OptProxies}

The tremendous advances in machine learning, especially deep learning, have enabled the modeling of highly complex relationships for accurate prediction. Recently, this has also spurred a growing interest in learning directly the mapping from the input parameters of an optimization problem to its optimal solution. In this setting, the required dataset is  $\{(\bar{u}_i, x^*_i)\}_{i \in \mathcal{I}}$, where $x^*_i$ denotes an optimal solution to the original problem \eqref{eqn:OriginalOP} for the input $\bar{u}_i$. The goal is to learn a model $m(\cdot) \in \mathcal{M}$ that returns an optimal (or near-optimal) solution for a given input by solving the following expected risk minimization problem:
\begin{equation}
\label{eqn:OptProxy}
    \minimize_{m(\cdot) \in \mathcal{M}} \; \frac{1}{N} \sum_{i \in \mathcal{I}} \ell \left( x^*_i, m(\bar{u}_i) \right).
\end{equation}
The corresponding surrogate model is simply
\begin{equation}
    x = m(u),
\end{equation}
which arrives at the solution through a function evaluation rather than by solving an optimization problem. Thus, $m$ serves as a computationally efficient proxy for the original optimization model. Importantly, as indicated by the loss function in problem \eqref{eqn:OptProxy}, this approach explicitly aims to minimize the decision prediction error. This is in contrast to the methods reviewed in Subsection \ref{sec:embedded}, which minimize the prediction error for individual constraint function values but cannot provide any guarantees in terms of how it affects the decision prediction error.

Although, in theory, any type of machine learning model (specified through $\mathcal{M}$) can be used in problem \eqref{eqn:OptProxy}, the preferred choice in the literature has been deep neural networks \citep{Sun2018, Pan2019, Krishnamoorthy2019, Karg2020, Kumar2021}. A major challenge in deep learning is that it is difficult to enforce constraints on the predictions, which in this case often leads to predicted solutions that are infeasible in the original optimization problem. \citet{Zamzam2020} address this challenge, when constructing neural network proxies for the AC optimal power flow (OPF) problem, by generating a training set of strictly feasible solutions through a modified AC OPF formulation and by using the natural bounds of the sigmoid activation function to enforce generation and voltage limits. Van Hentenryck and coworkers \citep{Fioretto2020} apply Lagrangian duality to consider constraints in their proposed deep learning framework, where the loss function in \eqref{eqn:OptProxy} is augmented with penalty terms derived from the Lagrangian dual of the original optimization problem and the corresponding Lagrange multipliers are updated using a subgradient method. This approach has been applied in several applications, including AC OPF \citep{Fioretto2020}, security-constrained OPF \citep{Velloso2021}, and job shop scheduling \citep{Kotary2022}. Another remedy is to correct an infeasible prediction by projecting it onto a suitably chosen set such that the projection is a feasible solution. For example, in MPC, such a set can be derived from the maximal control invariant set of the system to ensure recursive feasibility \citep{Chen2018a, Paulson2020}. 


While the development of the machine learning approaches described above is a more recent trend, exact methods for the construction of optimization proxies have been studied in the area of multiparametric programming for a long time. In multiparametric programming, optimal solutions are derived as explicit functions of the model parameters; these functions (or policies) can change depending on the region in which in the specific parameter vector lies. The goal is to determine these so-called critical regions and the optimal policy associated with each critical region. Once these are obtained offline, the online optimization reduces to simply selecting and applying the right function from a look-up table. This approach forms the basis for explicit MPC, which has been successfully applied in various real-world settings \citep{Alessio2009, Pistikopoulos2015}. However, a major challenge in multiparametric programming is the curse of dimensionality as the number of critical regions grows exponentially with the problem size. We refer the reader to \citet{Oberdieck2016} for a review of the extensive literature on multiparametric programming. 

\begin{remark}
In addition to the methods reviewed in this section, there is an extensive body of literature on the use of surrogate-based derivative-free optimization (DFO) techniques to solve complex optimization problems \citep{Rios2013}. We have consciously omitted referencing such work here. This is because our focus is specific, namely surrogate modeling for online optimization applications, where the same optimization model must be frequently solved with different values for its input parameters; hence, the surrogate model's parametric nature is crucial. The DFO approach lacks parametric characteristics such that one must solve the problem using the chosen DFO strategy every time the model parameters change, which means that the construction of the surrogate models within the DFO algorithm would have to happen online. As a result, DFO is typically not efficient enough to allow its use in real-time applications. On the other hand, parametric methods, like the ones reviewed here, offer a distinct advantage. With these methods, one trains a surrogate optimization model just once offline, and it can then be employed online to solve for any inputs $u \in \mathcal{U}$ without the need of retraining the surrogate model.
\end{remark}

\section{Decision-focused surrogate modeling}
\label{sec:DFSM}

The surrogate modeling frameworks reviewed in Section \ref{sec:Background} (also summarized in Figure \ref{fig:surrmod_frameworks}) all have their advantages and disadvantages. Optimization with embedded surrogate models is an intuitive approach that allows preservation of much of the structure of the original model. Domain knowledge can be effectively leveraged since for someone familiar with the physical system, it is usually easy to identify the part of the model that needs to be replaced as well as determine a suitable structure for the corresponding surrogate model. Often, only a small number of constraints are complicating; in that case, simply keeping the remaining constraints can help ensure feasibility. However, surrogate models generated using these methods may lead to solutions that are quite different from the optimal solutions to the original problem. This shortcoming can be overcome by learning optimization proxies in a decision-focused fashion. Using tailored deep learning architectures, this approach can achieve highly accurate and fast surrogate models. However, it is often less data-efficient, and cannot easily incorporate safety-critical hard constraints. In the following, we introduce the concept of decision-focused surrogate modeling (Figure \ref{fig:surrmod_frameworks}c), where we try to combine the desirable characteristics of the existing surrogate modeling approaches. We further present an inverse optimization approach to constructing such surrogate models for a certain class of problems.

\newcommand{\addpica}{\begin{minipage}{0.9\linewidth}
    \resizebox{\linewidth}{!}{\includegraphics{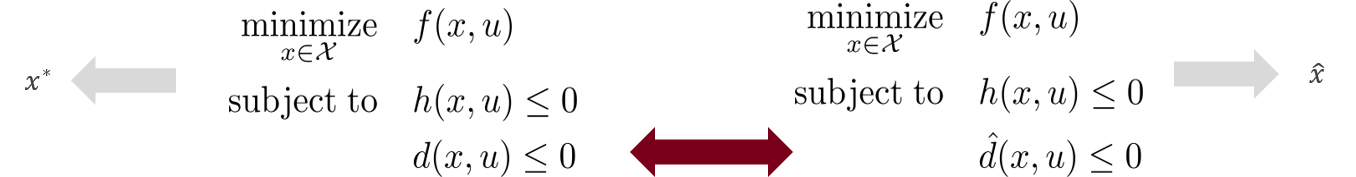}}%
    \caption*{(a) Optimization with embedded surrogate models}%
    \end{minipage}}
\newcommand{\addpicb}{\begin{minipage}{0.9\linewidth}
    \resizebox{\linewidth}{!}{\includegraphics{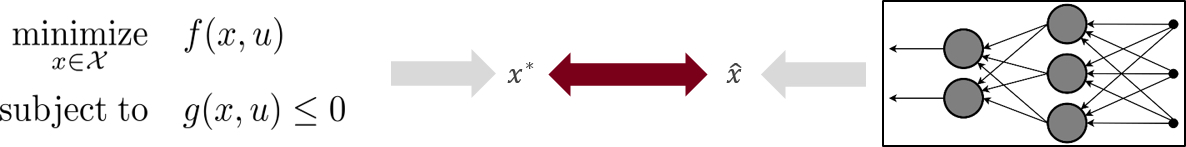}}%
    \caption*{(b) Learning optimization proxies}%
    \end{minipage}}
\newcommand{\addpicc}{\begin{minipage}{0.9\linewidth}
    \resizebox{\linewidth}{!}{\includegraphics{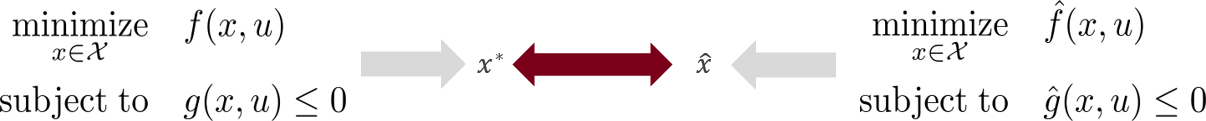}}%
    \caption*{(c) Decision-focused surrogate modeling}%
    \end{minipage}}
\begin{figure}
\begin{minipage}{\columnwidth}
        \centering
        \begin{tabular}{@{}c@{}}
           \toprule
            \addpica \\ \midrule
            \addpicb \\ \midrule
            \addpicc \\
            \bottomrule
        \end{tabular}
    \end{minipage}
    \caption{An overview of the surrogate modeling frameworks described in Section \ref{sec:Background} (a and b) and decision-focused surrogate modeling (c). The true optimal solution of the original optimization model is denoted by $x^*$, whereas $\hat{x}$ is the prediction generated by the surrogate. The gray arrows indicate the outputs of the various models while the red left-right arrows indicate what is considered in the loss function of each of the corresponding learning problems.}    \label{fig:surrmod_frameworks}
\end{figure}

\subsection{General formulation}

In the proposed framework, the data generation process is the same as in optimization proxy learning, where the dataset is $\{(\bar{u}_i, x^*_i)\}_{i \in \mathcal{I}}$ with $x^*_i$ denoting the true optimal solution to the original problem with input $\bar{u}_i$. Given such data, we directly train a surrogate optimization model defined by objective function $\hat{f}$ and constraint functions $\hat{g}$ that minimizes the decision prediction error. This learning problem can be formulated as follows:
\begin{varsubequations}{DFSLP}
\label{eqn:DFSM}
    \begin{align}
        \minimize_{\hat{f}(\cdot) \in \mathcal{F}, \, \hat{g}(\cdot) \in \mathcal{G}} \quad & \frac{1}{N} \sum_{i \in \mathcal{I}} \ell \left( x^*_i, \hat{x}_i \right) \\
        \st \quad \; & \hat{x}_i \in \argmin_{\tilde{x} \in \mathcal{X}} \{ \hat{f}(\tilde{x}, \bar{u}_i): \hat{g}(\tilde{x}, \bar{u}_i) \leq 0 \} \quad \forall \, i \in \mathcal{I}, \label{eqn:DFSM_argmin}
    \end{align}
\end{varsubequations}
where $\hat{f}(\cdot)$ and $\hat{g}(\cdot)$ can be chosen from some sets of functions $\mathcal{F}$ and $\mathcal{G}$, respectively. The solution predicted by the surrogate optimization model for input $\bar{u}_i$ is denoted by $\hat{x}_i$; hence, the loss function is analogous to the one in the optimization proxy learning problem \eqref{eqn:OptProxy}. As a result, this approach is decision-focused while allowing full flexibility in specifying the structure of the surrogate optimization model, including which original constraints to keep.

Note that since \eqref{eqn:DFSM} replaces the original constraints $g$ in \eqref{eqn:OriginalOP} with their surrogates $\hat{g}$, solving the resulting surrogate optimization problem could lead to solutions that are infeasible to \eqref{eqn:OriginalOP}. This presents a trade-off for enhancing computational efficiency. This trade-off is also present in many other methods that we reviewed in Section \ref{sec:Background}. Typical strategies for handling infeasible solutions involve projecting the predicted solution to the closest point on the feasible set or using the prediction to warm-start a fast local solver. With our surrogate modeling approach, infeasibility can be eased by retaining most of the original constraints and replacing only those posing computational challenges (similar to methods outlined in Section \ref{sec:embedded}). In a subsequent paper \citep{gupta2022decision}, we introduce a robust optimization method to identify regions of potential infeasibility within $\mathcal{U}$, which are then used to improve the surrogate model by acquiring additional training samples from these identified regions.


\begin{remark}
A closely related framework to our proposed decision-focused surrogate modeling approach is \textit{decision-focused learning} \citep{Wilder2019}. However, it is not motivated by the need for fast online optimization, but instead, it addresses the following problem: In traditional data-driven optimization, we often follow a two-stage predict-then-optimize approach, i.e., we first predict unknown input parameters $u$ from data with external features $r$ and then solve the optimization problem with those predicted inputs $u$. Here, the learning step focuses on minimizing the parameter estimation error; however, this does not necessarily lead to the best decisions (evaluated with the true parameter values) in the optimization step. In contrast, \textit{decision-focused learning}, also known as \textit{smart predict-then-optimize} \citep{Elmachtoub2022}, \textit{predict-and-optimize} \citep{Mandi2020}, and \textit{end-to-end learning for optimization} \citep{Donti2017}, integrates the two steps to explicitly account for the quality of the optimization solution in the learning of the model parameters.
\end{remark}


\subsection{An inverse optimization approach}

Problem \eqref{eqn:DFSM} can be viewed as a data-driven inverse optimization problem. Given an observed decision made by some agent, the goal of inverse optimization \citep{Ahuja2001} is to determine an optimization problem whose optimal solution matches and hence best explains the agent's decision; traditionally, only the objective function of the optimization model is assumed to be unknown \citep{ChanReviewPaper}. In the (noisy) data-driven setting, multiple decisions for different inputs are observed, and the goal is to find an optimization problem whose optimal solutions match the observations as closely as possible \citep{Gupta2021}. Applying this interpretation to the decision-focused surrogate modeling problem, the original optimization problem acts as the agent, and for each observation $i$, $x^*_i$ represents the observed decision, and $\bar{u}_i$ is the corresponding input; the surrogate optimization problem is then the optimization problem that we try to find such that its optimal solutions best resemble the observations.

Interpreting \eqref{eqn:DFSM} as an inverse optimization problem allows us to leverage existing methods from that literature to solve the problem. Recently, we developed an efficient data-driven inverse optimization framework that can incorporate both unknown objective functions and constraints \citep{Gupta2023}. This approach can be readily applied to a certain class of decision-focused surrogate modeling problems, which we formally define in the following. Here, we consider the original optimization problem to be a generally nonconvex nonlinear program (NLP) of the following form:
\begin{equation}
\label{eqn:OP}
\tag{OP}
\begin{aligned}
  \minimize_{x \in \mathbb{R}^n} \quad & f(x,u) \\
  \st \quad & g(x,u) \leq 0 \\
  & h(x,u) = 0.
\end{aligned}
\end{equation}

We assume that \eqref{eqn:OP} can be solved to generate the dataset $\mathcal{I}$ in \eqref{eqn:DFSM}. We use this dataset to learn a surrogate optimization problem which is a strictly convex optimization problem of the following form:

\begin{equation}
\label{eqn:SP}
\tag{SP}
\begin{aligned}
  \minimize_{x \in \mathbb{R}^n} \quad & \hat{f}(x,u;\theta) \\
  \st \quad & \hat{g}(x,u;\omega) \leq 0 \\
  & \hat{h}(x,u;\omega) = 0, 
\end{aligned}
\end{equation}
where functions $\hat{f}$ is strictly convex, $\hat{g}$ are convex, and $\hat{h}$ are linear in $x$. In \eqref{eqn:SP}, we parameterize the objective function and constraints with parameters $\theta$ and $\omega$, respectively. This parameterized surrogate problem can be learned in a decision-focused manner by solving \eqref{eqn:DFSM} with $\theta$ and $\omega$ as its decision variables; we use the term \textit{decision-focused surrogate optimization model} (DFSOM) to describe the resulting \eqref{eqn:SP}. In this work, we use mean squared error as the loss function, i.e. we solve \eqref{eqn:DFSM} with $\ell \left( x^*_i, \hat{x}_i \right) = \lVert x^*_i - \hat{x}_i \rVert_2^2$.

\subsection{Solution method}
The learning problem \eqref{eqn:DFSM} is a bilevel program \citep{dempe2020bilevel} with multiple lower-level optimization problems. Despite the convexity of the lower-level problems, \eqref{eqn:DFSM} is difficult to solve due to its poor scalability \citep{Aswani2018}. In our previous work \citep{Gupta2023}, we addressed this problem with an efficient decomposition-based strategy that generates high-quality solutions; the same algorithm is employed here. For the sake of brevity, we only provide a brief overview of our solution method in this subsection and refer the interested reader to Gupta and Zhang \cite{Gupta2023} for a more detailed discussion.

Our solution approach involves reformulating the bilevel problem  into a single-level optimization problem by replacing the lower-level problems with their \emph{Karush-Kuhn-Tucker} (KKT) optimality conditions. It is important to highlight that this reformulation requires that \eqref{eqn:SP} satisfies some constraint qualification conditions. This requirement can be met by designing \eqref{eqn:SP} with all linear constraints (as done in our computational experiments in Section \ref{sec:case_studies}), or by ensuring that the set defined by $\hat{g} < 0$ and $\hat{h} = 0$ is always feasible. The reformulation of \eqref{eqn:DFSM} using KKT conditions results in an NLP as follows:

\begin{subequations}
    \label{eqn:KKT-Convex}
        \begin{align}
        \minimize\limits_{\theta \in \Theta, \, \omega \in \Omega, \hat{x} ,\lambda, \mu} \quad & \sum\limits_{i \in \mathcal{I}} \, \lVert x^*_i - \hat{x}_i \rVert_2^2 \\
        \st \quad \; \; &  \nabla {\hat{f}} (\hat{x}_i, \bar{u}_i; {\theta}) + \lambda_i^{\top} \nabla {\hat{g}} (\hat{x}_i, \bar{u}_i; {\omega}) + \mu_i^{\top} \nabla {\hat{h}} (\hat{x}_i, \bar{u}_i; {\omega}) = 0 \quad \forall \, i \in \mathcal{I} \label{eqn:stat} \\
        & {\hat{g}}(\hat{x}_i, \bar{u}_i; {\omega}) \leq 0 \quad \forall \, i \in \mathcal{I} \label{eqn:PF1} \\
        & {\hat{h}}(\hat{x}_i, \bar{u}_i; {\omega}) = 0 \quad \forall \, i \in \mathcal{I} \label{eqn:PF2} \\
        & \lambda_{i}^{\top} \, {\hat{g}}(\hat{x}_i, \bar{u}_i;{\omega}) = 0 \quad \forall \, i \in \mathcal{I}  \label{eqn:CS} \\
        & \lambda_i \geq 0, \hat{x}_i \in \mathbb{R}^n \quad \forall \, i \in \mathcal{I},
        \end{align}
    \end{subequations}
where the dual variables for the inequality and equality constraints of the lower-level problems in \eqref{eqn:DFSM} are respectively denoted by $\lambda$ and $\mu$. Constraints \eqref{eqn:stat}, \eqref{eqn:PF1}-\eqref{eqn:PF2}, and \eqref{eqn:CS} formulate the stationarity, primal feasibility, and complementary slackness conditions, respectively, for the lower-level problems. 

The reformulated problem is a nonconvex NLP with complementarity constraints that violate standard constraint qualification conditions. As these problems are known to cause convergence difficulties for NLP solvers, we further consider an \emph{exact penalty reformulation} \citep{nocedal2006numerical} of \eqref{eqn:KKT-Convex}. This reformulation yields the following problem with a regularized feasible region that satisfies the necessary constraint qualifications if the description of sets $\Theta$ and $\Omega$ also satisfy necessary regularity conditions: 
\begin{equation}
    \label{eqn:KKT-EPR}
        \begin{aligned}
        \minimize\limits_{{\theta} \in \Theta, \, {\omega} \in \Omega, \hat{x} ,\lambda, \mu} \quad & \, \sum_{i \in \mathcal{I}}  \lVert x^*_i - \hat{x}_i \rVert_2^2 
        \\ & + c^{\top} 
        \underbrace{\begin{bmatrix} 
        \sum\limits_{i \in \mathcal{I}} \lvert \nabla \hat{f} (\hat{x}_i, \bar{u}_i; {\theta}) + \lambda_i^{\top} \nabla {\hat{g}} (\hat{x}_i, \bar{u}_i; {\omega}) + \mu_i^{\top} \nabla {\hat{h}} (\hat{x}_i, \bar{u}_i; {\omega}) \rvert \\
        \sum\limits_{i \in \mathcal{I}} \mathrm{max} \{ 0, \hat{g}(\hat{x}_i, \bar{u}_i; {\omega}) \} \\
        \sum\limits_{i \in \mathcal{I}} \lvert \hat{h}(\hat{x}_i, \bar{u}_i; {\omega}) \rvert \\
        \sum\limits_{i \in \mathcal{I}} \lvert \lambda_{i}^{\top} \, \hat{g}(\hat{x}_i, \bar{u}_i; {\omega}) \rvert
        \end{bmatrix}}_{P} \\
        \st \quad & \lambda_i \geq 0, \hat{x}_i \in \mathbb{R}^n \quad \forall \, i \in \mathcal{I},
        \end{aligned}
    \end{equation}
where $c$ are the positive penalty parameters. The exactness of this reformulation holds only for penalty parameters greater than a certain threshold, which is hard to determine a priori in practice. As a result, we take an iterative approach, starting with small values for $c$ and gradually increasing them in subsequent iterations (by a factor of $\rho$) until the solution of \eqref{eqn:KKT-EPR} is also a feasible solution for \eqref{eqn:KKT-Convex}. We use a feasibility tolerance of $\epsilon$ as the termination criterion for the algorithm. A pseudocode that summarizes our overall approach is shown in Algorithm \ref{alg:PBCD}.

While \eqref{eqn:KKT-Convex} can be solved using the penalty reformulation, the reformulated problem \eqref{eqn:KKT-EPR} remains a nonconvex NLP whose large instances are generally difficult to solve. Here, we notice that for certain classes of \eqref{eqn:SP}, including quadratic programs (QPs), \eqref{eqn:KKT-EPR} becomes a multiconvex optimization problem \citep{Shen2017} (MCP) if the sets $\Theta$ and $\Omega$ are convex. This feature of \eqref{eqn:KKT-EPR} allows us to solve it with an efficient block-coordinate-descent (BCD) algorithm. We specifically highlight QPs as an example here because QPs are commonly used to model/approximate RTO and MPC problems. In fact, we show several examples of real-world systems in Section \ref{sec:case_studies} where we use a QP to approximate their original nonconvex RTO problem. Therefore, the applicability to QPs makes our decomposition scheme especially appealing for decision-focused surrogate modeling. 

\begin{remark}
The objective function of \eqref{eqn:KKT-EPR} contains nonsmooth terms due to the use of $\ell_1$-norm-based penalty functions. This problem can be exactly reformulated into one with a smooth objective function by introducing additional variables to linearize the penalty terms. The reformulation results in an increased size of the problem; however, in our previous work\citep{Gupta2023}, we show that \eqref{eqn:KKT-EPR} provides significantly better solutions than when IPOPT is applied directly on \eqref{eqn:KKT-Convex}.  
\end{remark}

\begin{algorithm}
    \begin{algorithmic}[1]
        \State initialize: $ k \gets 1, ({\theta}, {\omega}, \hat{x}, \lambda, \mu) \gets ({\theta}_0, {\omega}_0, \hat{x}_0, \lambda_0, \mu_0)$ and $c \gets c_1$
        \While{$\lVert P \rVert > \epsilon $}
            \State solve \eqref{eqn:KKT-EPR} with BCD or an NLP solver (e.g. IPOPT \citep{wachter2006})
            \State $c_{k+1} \gets c_k + \rho c_k$
            \State $k \gets k+1$
        \EndWhile
        \State \Return $({\theta}_k, {\omega}_k, \hat{x}_k, \lambda_k, \mu_k)$
        \caption{A penalty block coordinate descent algorithm for solving \eqref{eqn:KKT-Convex}.}
        \label{alg:PBCD}
    \end{algorithmic}
\end{algorithm}

\section{Computational case studies} \label{sec:case_studies}
We present numerical results from three case studies based on typical nonlinear chemical engineering systems where we apply the proposed decision-focused surrogate modeling strategy to simplify the optimization task. The first two case studies are based on single-input systems for which we present detailed analyses demonstrating how our approach can provide excellent decision prediction accuracy with a much simpler model than is typically required by traditional methods that construct surrogate optimization problems with embedded surrogate models. The third case study is based on a larger multi-input multi-output blending network system where we show the data efficiency of our approach in comparison to black-box optimization proxies. We also demonstrate the effectiveness of the resulting DFSOM in reducing the computational burden of optimizing the system in real time. All computer code along with the datasets used for these case studies is available at https://github.com/ddolab/DecFocSurrMod.

\subsection{Real-time optimization of a continuous stirred tank reactor (CSTR) system}
We consider the problem of real-time optimization of a CSTR operating in a stochastic environment. Specifically, we optimize the operation of an ideal adiabatic CSTR \citep{economou1986} that is carrying out an exothermic reversible reaction between reactant $A$ and product $R$. The concentration of the inlet stream, which does not contain any $R$, is subject to observable disturbances. Here, the primary goal of RTO is to maximize the concentration of the product $R$ in the outlet stream by manipulating the temperature of the inlet stream. This can be formulated as the following optimization problem:
\begin{equation}
    \label{eqn:cstr_rto}
    \begin{aligned}
        \maximize_{T_i, T_o, A_o, R_o} \quad & 2.009 R_o - 1.657 \times 10^{-3} (T_i - 410)^2 \\
        \st \quad & 0 = 1/\tau (A_i - A_o) - k_1 (T_o) A_o + k_{-1}(T_o) R_o \\ 
        & 0 = - R_o/\tau + k_1(T_o) A_o - k_{-1}(T_o) R_o \\
        & 0 = 1/\tau (T_i - T_o) + \frac{-\Delta H_R}{ \rho C_p} \left(k_1(T_o) A_o - k_{-1}(T_o)R_o \right),
    \end{aligned}
\end{equation}
where $A_i$ is the concentration of $A$ in the inlet stream, and $T_i$ is the inlet temperature. Similarly, variables with the subscript $o$ characterize the properties of the outlet stream. In \eqref{eqn:cstr_rto}, the first two constraints are the mass balances whereas the third constraint specifies the energy balance for the reactor. The paper by Economou et al. \citep{economou1986} provides details of all model parameters used in this case study.

Problem \eqref{eqn:cstr_rto} is a nonconvex NLP due to the dependence of the forward and backward reaction rate constants ($k_1$ and $k_{-1}$) on the reactor temperature ($T_o$). We use decision-focused surrogate modeling to learn a convex surrogate optimization model for the RTO problem. We achieve this by replacing the original rate law expression by the following approximation that is linear in the decision variables of \eqref{eqn:cstr_rto}:
\begin{equation} \label{eqn:ratelaw_surr}
k_1(T_o) A_o - k_{-1}(T_o) R_o \; \rightarrow \; a(A_i) \, T_o + b(A_i) \, A_o + c(A_i) \, R_o, 
\end{equation}
where $a$, $b$, and $c$ are functions of the input parameter $A_i$. In addition, we also learn scalar coefficient values for a new quadratic objective function for the DFSOM.

\subsubsection{Surrogate model complexity and data efficiency}
\label{sec:CSTRsub}
We begin by examining how using different functional forms for the coefficients $a(A_i)$, $b(A_i)$, and $c(A_i)$ in \eqref{eqn:ratelaw_surr} affects the accuracy of the learned DFSOM. To simplify our analysis, we focus on these coefficients being polynomial functions of the input disturbance value. We test our model's accuracy by varying the complexity of the functions from constants up to cubic polynomials of $A_i$. In each case, we train the DFSOM using a dataset $\mathcal{I}$ of 1,000 samples, with each sample containing the input parameter $A_i$ and the corresponding optimal $(T_i, A_o, R_o)$ values obtained by solving the original NLP \eqref{eqn:cstr_rto}. After training, we evaluate the accuracy of the model by testing it on a separate dataset of 100 samples. 

The results of the model complexity analysis are shown in Figure \ref{fig:model_complexity}. In the plot, the "Sparse 3" column represents the case in which we specify the coefficient models to be cubic but add a term to the objective function of \eqref{eqn:DFSM} that penalizes the sum of absolute values of the polynomial coefficients. This technique is commonly used in machine learning to control model complexity. We find that this approach results in simple yet highly accurate polynomials, so we choose "Sparse 3" as the model for our further analysis.

Figure \ref{fig:model_complexity} shows that increasing the model complexity generally leads to a more accurate optimization problem. However, it is important to note that this added complexity lies in the input parameter space, and the optimization problem in the decision variable space stays convex. This is one of the unique features of our approach, allowing us to transfer the complexity from the decision variable space to the input space while still learning the relevant characteristics of the original optimization model. None of the traditional methods that employ embedded surrogate models can learn a similar model because their model learning step is entirely separate from the optimization step, so there is no differentiation between input parameter space and decision variable space when constructing the surrogate model.

We then investigate how the size of the training dataset affects the accuracy of our model. This is important because constructing the dataset requires solving a complex optimization problem. Figure \ref{fig:train_data} shows that our method can produce an accurate and robust model with as few as 25 samples, indicating that a large dataset may not be necessary. In the next subsection, we compare the performance and data efficiency of our method with traditional methods using embedded surrogate models. 

\usepgfplotslibrary{colorbrewer}
\pgfplotsset{compat = 1.16, cycle list/Set1-8} 
\usetikzlibrary{pgfplots.statistics, pgfplots.colorbrewer} 

\begin{figure}[h!]
    \centering
    \subfloat[Effect of model complexity]{
        \label{fig:model_complexity}
        \begin{tikzpicture}[scale=0.9]
            \pgfplotstableread[col sep=comma]{Figures/cstr_model_comp.csv}{\csvdata}
            \begin{axis}[
                boxplot/draw direction = y,
                x axis line style = {opacity=100, thick},
                y axis line style = {thick},
                enlarge y limits,
                ymajorgrids,
                xtick = {1, 2, 3, 4, 5},
                xticklabel style = {align=center},
                xticklabels = {$0$,$1$,$2$,$3$,Sparse $3$},
                xtick style = {}, 
                ylabel = {mean $ \lvert (T_i^*) - (\hat{T_i}) \rvert $},
                xlabel = {Degree of polynomial}
            ]
                \foreach \n in {0,...,4} {
                    \addplot+[boxplot = {every median/.style={red,ultra thick}}, draw=blue] table[y index=\n] {\csvdata};
                }
            \end{axis}
        \end{tikzpicture}}
    \subfloat[Effect of training dataset size]{
        \label{fig:train_data}
        \begin{tikzpicture}[scale=0.9]
            \pgfplotstableread[col sep=comma]{Figures/full_pred_error.csv}{\csvdata}
            \begin{axis}[
                boxplot/draw direction = y,
                x axis line style = {opacity=100, thick},
                y axis line style = {thick},
                enlarge y limits,
                ymajorgrids,
                xtick = {1, 2, 3, 4, 5, 6, 7, 8},
                xticklabel style = {align=center},
                xticklabels = {$5$, $10$, $25$, $50$, $75$, $100$, $250$, $500$, $1000$},
                xtick style = {}, 
                ylabel = {mean $ \lvert (T_i^*) - (\hat{T_i}) \rvert $},
                ymax = 1.50,
                xlabel = {$|\mathcal{I}|$}
            ]
                \foreach \n in {0,...,7} {
                    \addplot+[boxplot = {every median/.style={red,ultra thick}}, draw=blue] table[y index=\n] {\csvdata};
                }
            \end{axis}
        \end{tikzpicture}}
    \caption{\label{fig:cstr_CustAdaptive} The true optimal solution obtained by solving \eqref{eqn:cstr_rto} is denoted by $T_i^*$, and $\hat{T}_i$ represents the solution of the DFSOM. The box plots show the interquartile ranges of prediction error for the surrogate models obtained with ten different instances of randomly generated training data. We use a separate test dataset of 100 samples to compute the prediction error.}
    \end{figure}
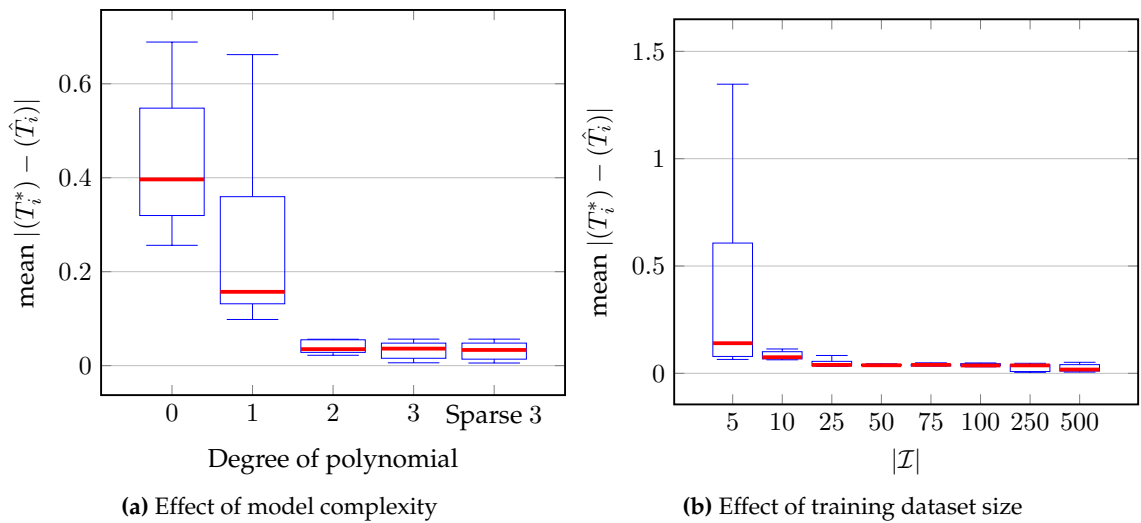

\subsubsection{Comparison with traditional embedded surrogate model approach}

To compare our proposed strategy to other surrogate modeling methods popular in PSE, we construct surrogate models for the rate law expression. We use a large training dataset of 3,000 points, where each point consists of $(A_o, R_o, T_o)$ as inputs and the true rate value as the output. These surrogate models learn a function $\hat{r}(A_o, R_o, T_o)$ that produces approximately the same reaction rate as the original model (on the left-hand side) in \eqref{eqn:ratelaw_surr}. We then substitute the original rate law in \eqref{eqn:cstr_rto} with $\hat{r}$ to obtain a surrogate optimization model, which we solve using the local NLP solver IPOPT \citep{wachter2006}. Note that although the resulting surrogate optimization models are generally nonconvex, we use a local solver because it is typically preferred in RTO settings due to the high computational effort required by global solvers.

We consider the following two surrogate modeling methods for the rate law expression:
\begin{itemize}
\item ES-ALAMO - We use ALAMO \citep{Cozad2014} (version 2022.10.7) to estimate $\hat{r}$. We allow all available basis functions except for the sine and cosine functions. We manually substituted the algebraic $\hat{r}$ expressions obtained using ALAMO into an RTO problem implementation in the JuMP \citep{dunning2017} modeling environment available in the Julia programming language. 

\item ES-NN - We train a neural network (NN) to estimate the reaction rate value. The NN model comprises three inputs, four hidden layers with ten nodes, and a single output node. We use smooth sigmoid activation functions for the hidden layers to ensure that the resulting optimization model remains solvable with IPOPT. We used the Python library TensorFlow \citep{tensorflow2015} v2.10.0 to train the NN models. These models were incorporated into the optimization problems using the OMLT \citep{ceccon2022} package v1.0.      
\end{itemize}

\begin{figure}[h!]
    \includegraphics[width=0.6\textwidth]{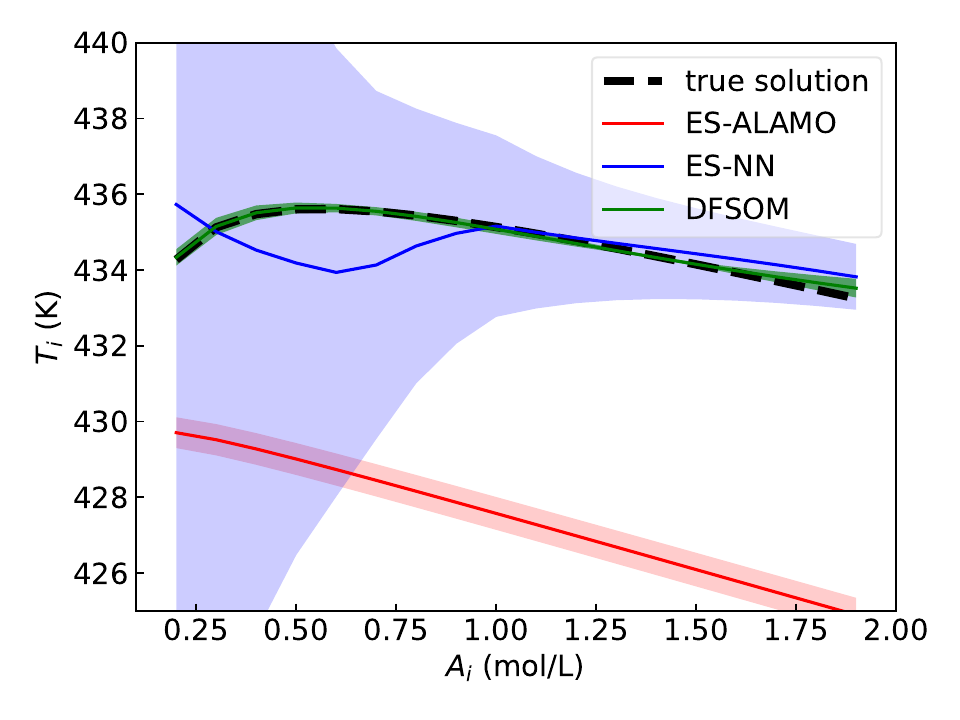}
    \centering
    \caption{Line plots show the mean optimal $T_i$ values obtained using different surrogate modeling methods, with shaded areas indicating one standard deviation. Five surrogate optimization models were developed for each method using different datasets.}
    \label{fig:cstr_comparison}
\end{figure}

Figure \ref{fig:cstr_comparison} shows a comparison of optimal $T_i$ values obtained using different methods. The DFSOM, trained with only 50 data points, significantly outperforms the two traditional methods. While ES-ALAMO yields a robust model that does not vary much with changing data, the predicted optimal $T_i$ values differ significantly from their true values. In contrast, while ES-NN yields mean $T_i$ values close to the true optimum, its output is highly sensitive to the training dataset used. The decision-focused model provides both high accuracy and robustness, independent of the quality of the dataset.
 
\begin{figure}[h!]
    \centering
    \begin{minipage}{\linewidth}
    \resizebox{\linewidth}{!}{\includegraphics{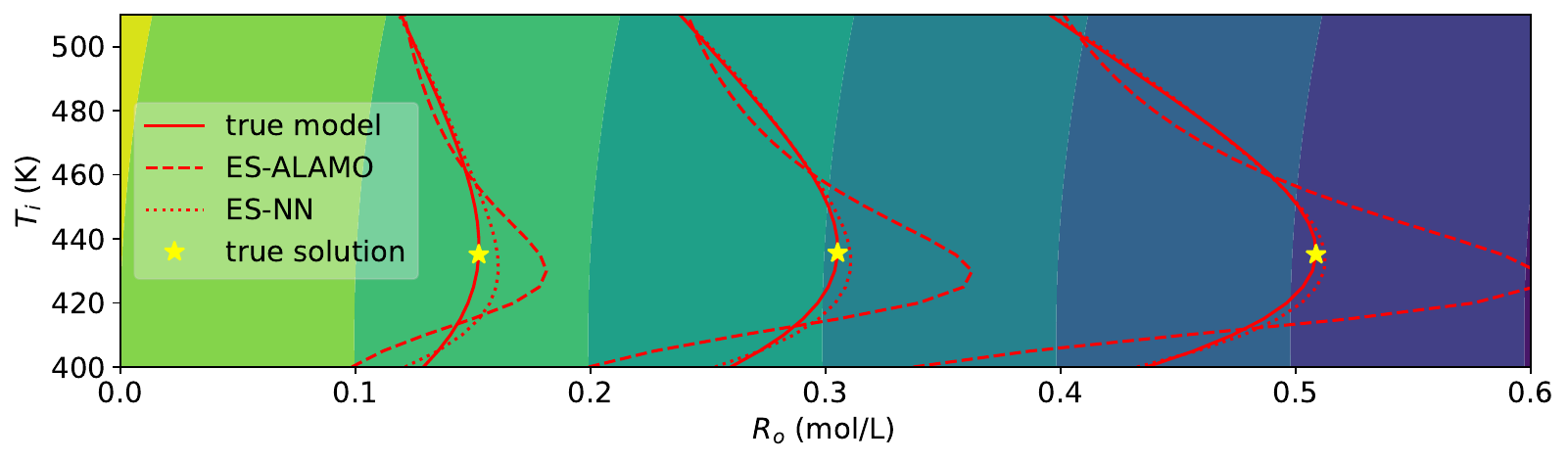}}%
    \caption*{(a)}%
    \end{minipage}%
    \hfill
    \begin{minipage}{\linewidth}
    \resizebox{\linewidth}{!}{\includegraphics{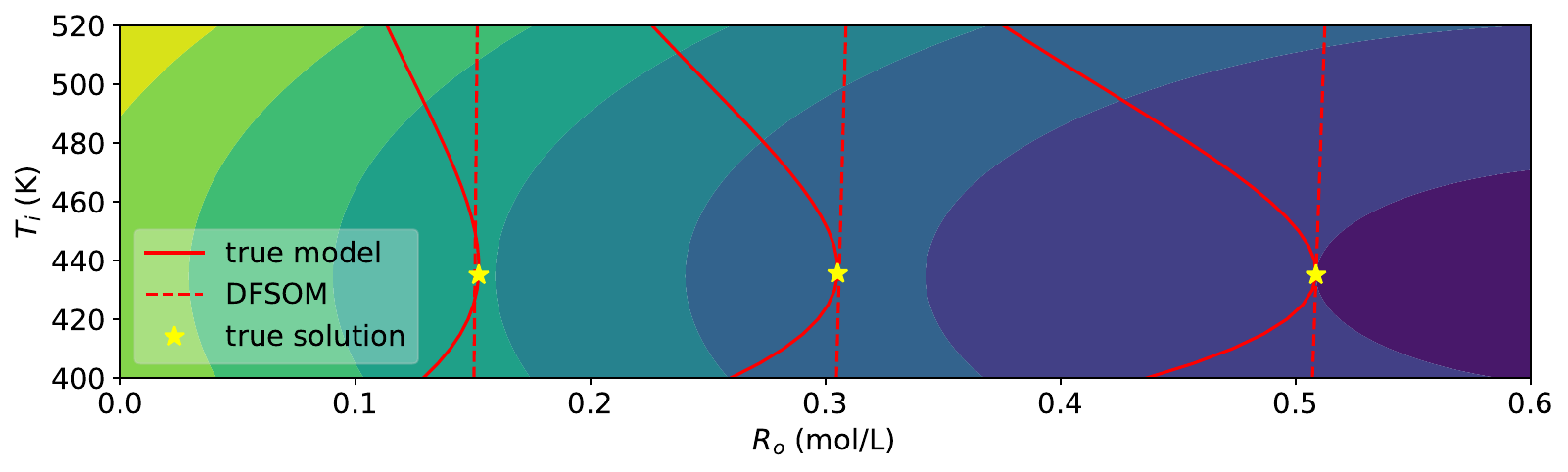}}%
    \caption*{(b)}%
    \end{minipage}%
    \caption{The feasible regions of the original and surrogate RTO problems projected onto the two-dimensional variable space between $R_o$ and $T_i$. From left to right, the three line plots depict the feasible regions for $A_i = 0.3, 0.6,$ and $1.0$ mol/L, respectively. The contours in the background represent the objective function for the true model in (a) and for the DFSOM in (b), with darker colors indicating larger values.}
    \label{fig:cstr_feasible_region}
\end{figure}

We further evaluate the quality of different models by comparing the feasible regions of the surrogate optimization models with the true feasible region. This analysis is shown in Figure \ref{fig:cstr_feasible_region}. Methods that embed surrogate models in the optimization problem, as compared in Figure \ref{fig:cstr_feasible_region}a, aim to exactly replicate the feasible region, often requiring complex models that need larger training datasets. However, even with a good surrogate model that closely replicates the true feasible region, the decision prediction error can still be high due to small discrepancies, as we observe for ES-NN. In contrast, decision-focused learning achieves high accuracy with simple models by focusing only on predicting the optimal solutions. As shown in Figure \ref{fig:cstr_feasible_region}b, although the feasible region of the DFSOM is different from the true feasible region, they coincide exactly at the optimum for the original model. The learning problem \eqref{eqn:DFSM} adjusts the surrogate optimization model's objective function to ensure this point is the optimal solution. Overall, we find that decision-focused learning allows us to significantly reduce the complexity of the surrogate optimization model while retaining its accuracy in predicting the optimal solutions.

\subsection{Real-time optimization of a heat exchanger network}

\begin{figure}[h!]
    \includegraphics[width=0.50\textwidth]{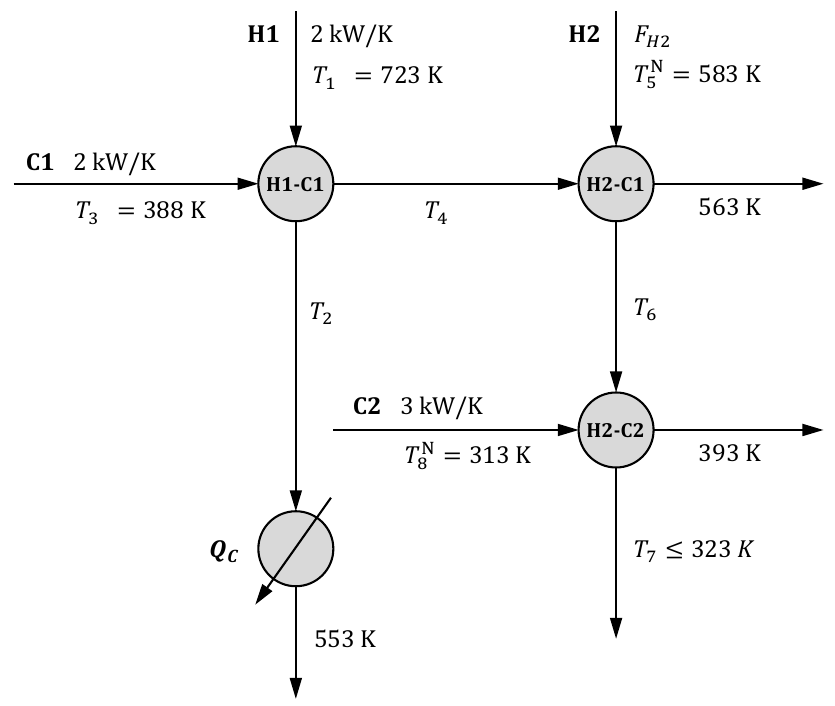}
    \centering
    \caption{Given heat exchanger network.}
    \label{fig:HXNet}
\end{figure}

Our second case study considers a heat exchanger network adapted from Biegler et al. \citep{biegler1997}, as shown in Figure \ref{fig:HXNet}. The inlet temperature of stream H2, denoted by $T_5$, has a nominal value of $583$ K but is subject to random disturbances. Upon a change in $T_5$, we optimize the operation of the heat exchanger network by solving the following nonconvex NLP:

\begin{subequations}
    \label{eqn:HXNet}
    \begin{align}
        \minimize\limits_{Q_c, \, F_{H2}} \quad & 10^{-2} \, Q_c + 4 \, (F_{H2} - 1.7)^2 \\
        \st \quad & 0.5 \, Q_c + 165 \geq 0 \\
        & -10 - Q_c + \left(T_5 - 558 + 0.5 \, Q_c \right) \, F_{H2} \geq 0 \label{eqn:HX_nonlinearC} \\
        & -10 - Q_c + \left( T_5 - 393 \right) \, F_{H2} \geq 0 \\
        & -250 - Q_c + \left( T_5 - 313 \right) \, F_{H2} \geq 0 \\
        & -250 - Q_c + \left( T_5 - 323 \right) \, F_{H2} \leq 0 \\
        & Q_c \geq 0, \, F_{H2} \geq 0,
    \end{align}
\end{subequations}
where the cooling duty $Q_c$ and the heat capacity flowrate $F_{H2}$ are adjustable variables. 

The nonconvexity in \eqref{eqn:HXNet} arises from the bilinear term in constraint \eqref{eqn:HX_nonlinearC}. Therefore, we apply our approach to replace the nonconvex term with the following approximation that is linear in the decision variables of \eqref{eqn:HXNet}:
\begin{equation}
    \label{eqn:LinearApp}
    Q_c \, F_{H2} \; \rightarrow \; a(T_5) \, Q_c + b(T_5) \, F_{H2},
\end{equation}
where $a$ and $b$ are functions of the input parameter $T_5$. Along with this, we also learn new coefficients for the quadratic objective function. Unlike in the previous case study, here these coefficients are also functions of the input parameter $T_5$. We specify the model for all unknown parameters as ``Sparse 3'' polynomials in $T_5$ (which we defined in Section \ref{sec:CSTRsub}). 

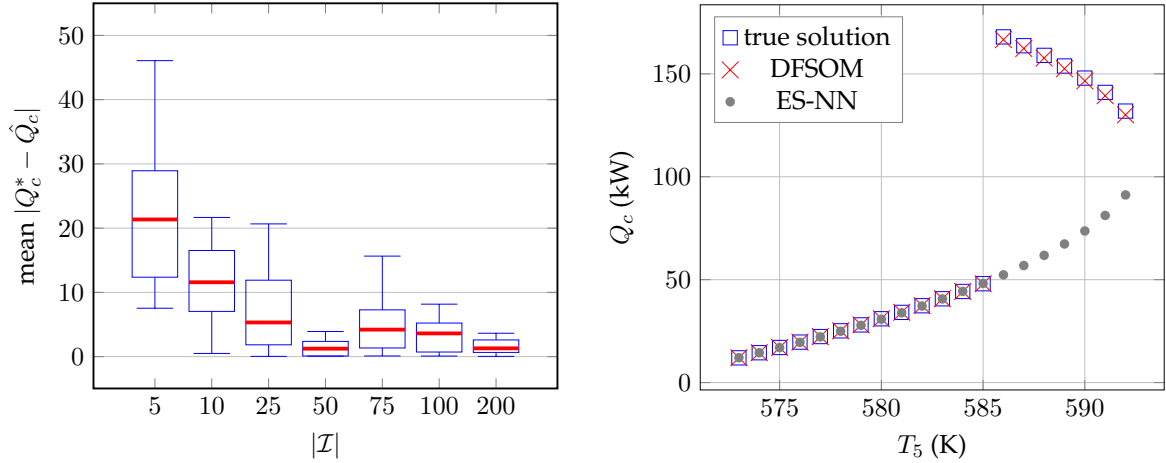
\begin{figure}[h!]
    \centering
    \subfloat[Effect of training dataset size on model accuracy. We use a separate test dataset of 100 samples to compute the prediction error.]{
        \label{fig:hx_data_eff}
        \begin{tikzpicture}[scale=0.9]
            \pgfplotstableread[col sep=comma]{Figures/hx_full_pred_error.csv}{\csvdata}
            \begin{axis}[
                boxplot/draw direction = y,
                x axis line style = {opacity=100, thick},
                y axis line style = {thick},
                enlarge y limits,
                ymajorgrids,
                xtick = {1, 2, 3, 4, 5, 6, 7},
                xticklabel style = {align=center},
                xticklabels = {$5$, $10$, $25$, $50$, $75$, $100$, $200$},
                xtick style = {}, 
                ylabel = {mean $\lvert Q_c^* - \hat{Q_c} \rvert $},
                ymax = 50,
                xlabel = {$|\mathcal{I}|$}
            ]
                \foreach \n in {0,...,6} {
                    \addplot+[boxplot = {every median/.style={red,ultra thick}}, draw=blue] table[y index=\n] {\csvdata};
                }
            \end{axis}
        \end{tikzpicture}} \hfil
    \subfloat[Comparison of decision-focused surrogate modeling with embedded surrogate models]{
        \label{fig:hxnet_comparison}
        \begin{tikzpicture}[scale=0.9]
            \usetikzlibrary[plotmarks]
            \begin{axis}[
                xlabel = {$T_5$ (K)},
                ylabel = {$Q_c$ (kW)},
                legend entries = {true solution, DFSOM, ES-NN},
                legend pos = north west,
                legend style = {draw=gray},
                grid=major
            ]
                \addplot[only marks, blue, mark size=3, mark=square] table [x=T5, y=true, col sep = comma] {Figures/hx_data.csv};
                \addplot[only marks, red, mark=x, mark size=5] table [x=T5, y=io, col sep = comma] {Figures/hx_data.csv};
                \addplot[only marks, gray, draw opacity=0, mark size=2] table [x=T5, y=nn, col sep = comma] {Figures/hx_data.csv};
            \end{axis}
        \end{tikzpicture}}
    \caption{\label{fig:hx_CustAdaptive} The true optimal solution obtained by solving \eqref{eqn:HXNet} is denoted by $Q_c^*$, and $\hat{Q}_c$ represents the solution of the surrogate optimization model. The box plots show the interquartile ranges of prediction error for the surrogate optimization models obtained with ten different instances of randomly generated training data.}
\end{figure}

To train the models for $a$, $b$, and the objective coefficients, we solve \eqref{eqn:DFSM} with datasets of varying sizes. Each data point in the dataset consists of an input $T_5$ and the corresponding optimal $Q_c^*$ and $F_{H2}^*$ values that we obtain by solving \eqref{eqn:HXNet}. We evaluate the performance of the resulting surrogate optimization models using a test dataset of 100 data points. The analysis of these results is available in Figure \ref{fig:hx_data_eff}, where we show the evolution of prediction quality as a function of the training dataset size. We find that the prediction error of the models converges to its minimum value with as few as 50 data points. Therefore, we set $|\mathcal{I}|$ to $50$ for our further analysis.

Next, we compare the performance of our decision-focused approach with the traditional approach with embedded surrogate models. To do this, we train an NN with two hidden layers, each consisting of ten nodes. The objective is to output the value of the bilinear term, given $Q_c$ and $F_{H2}$ as inputs. The dataset used for this purpose consists of 1,000 data points. Similar to the first case study, we maintain the sigmoid activation function for the hidden layers to guarantee the solvability of the surrogate optimization model with IPOPT. We then embed this NN in \eqref{eqn:HXNet} by replacing the bilinear term, resulting in a surrogate optimization model that we call ES-NN.

Due to the simplicity of the function that the NN intends to replace, the NN learns to approximate the bilinear term almost perfectly. However, as we show in Figure \ref{fig:hxnet_comparison}, ES-NN does not always yield the true optimal $Q_c$ values as its solution. This happens because the NN-based problem is nonconvex and IPOPT recovers one of the two local solutions that are possible for a given $T_5$ value. Here, a global solver can be used to remedy this situation, but this may not be practical for RTO due to the associated computational expense. In contrast, the DFSOM, which is convex by design, correctly identifies the existence of the discontinuity in the global optimal solution space with IPOPT. Additionally, we find that the DFSOM produces solutions that are nearly identical to those produced by the true problem for all values of $T_5$ considered in this case study.

In conclusion, this case study showcases a notable benefit of decision-focused surrogate modeling, whereby the surrogate optimization models resulting from this approach exhibit convexity, rendering them amenable to efficient local solvers. This obviates concerns regarding the convergence to suboptimal solutions, which is a common issue when using nonconvex models in RTO applications.

\subsection{Real-time optimization of a blending network}
We consider a network of blending nodes that mix material streams of various specifications to produce products with desired qualities. An optimization problem to minimize the cost of operation of this network is as follows \cite{adhya1999lagrangian}:   
\begin{subequations}
\label{eqn:blending_prob}
\begin{align}
    \minimize\limits_{f, x, q} \quad & \sum\limits_{j \in \mathcal{J}} \sum\limits_{i \in \mathcal{N}_j} c_{ij} \left( f_{ij} + \frac{1}{1000} (f_{ij} - \bar{f}_{ij})^2 \right) - \sum\limits_{k \in \mathcal{K}} d_{k} \sum\limits_{j \in \mathcal{J}} x_{jk} \label{eqn:blend_obj} \\
    \st \quad & \sum\limits_{i \in \mathcal{N}_j} f_{ij} = \sum\limits_{k \in \mathcal{K}} x_{jk} \quad \forall \, j \in \mathcal{J} \label{eqn:node_mb} \\
    & q_{jw} \sum\limits_{k \in \mathcal{K}} x_{jk} = \sum\limits_{i \in \mathcal{N}_j} \lambda_{ijw} f_{ij} \quad \forall \, j \in \mathcal{J}, w \in \mathcal{W} \label{eqn:node_quality} \\
    & \sum\limits_{j \in \mathcal{J}} x_{jk} \leq S_k \quad \forall \, k \in \mathcal{K} \label{eqn:demand_satis} \\
    & \sum\limits_{j \in \mathcal{J}} q_{jw} x_{jk} \leq Z_{kw} \sum\limits_{j \in \mathcal{J}} x_{jk} \quad \forall \, k \in \mathcal{K}, w \in \mathcal{W} \label{eqn:quality_satis} \\
    & f_{ij} \geq 0 \quad \forall \, j \in \mathcal{J}, i \in \mathcal{N}_j \\
    & q_{jw} \geq 0 \quad \forall \, j \in \mathcal{J}, w \in \mathcal{W} \\
    & x_{jk} \geq 0 \quad \forall \, j \in \mathcal{J}, k \in \mathcal{K},
\end{align}
\end{subequations}
where $\mathcal{J}$ is the set of blending nodes, and $\mathcal{N}_j$ are the sets of material streams entering node $j$; we use variable $f_{ij}$ to denote the flow rate of stream $i$ entering node $j$. The product streams are generated by combining the outflows from different blending nodes and are represented by the set $\mathcal{K}$. The outflow rate from node $j$ to product stream $k$ is denoted by the variable $x_{jk}.$ Additionally, the set $\mathcal{W}$ contains the various components present in the material streams whose specifications need to be maintained in the output streams. We determine the quality of node $j$ for a specific component $w$ using the variable $q_{jw}$.

In problem \eqref{eqn:blending_prob}, we enforce mass balance for each node through constraints \eqref{eqn:node_mb}. The quality of a blending node $j$ is determined as a function of the specifications, $\lambda_{ijw}$, of the streams entering that node, through constraints \eqref{eqn:node_quality}. Additionally, we ensure that the outflow rate of a product stream $k$ does not exceed its market demand, $S_k$, using constraints \eqref{eqn:demand_satis}. Finally, constraints \eqref{eqn:quality_satis} set the specification of component $w$ in product stream $k$ below its permissible level, $Z_{kw}$. The objective function of \eqref{eqn:blending_prob} comprises of three terms. The first represents the linear cost of acquiring a material stream, while the second penalizes the deviation of inlet flow rates from their nominal values, $\bar{f}_{ij}$. The third term represents the revenue generated through the product streams.

We consider a scenario in which the RTO of the blending network operation is desired in the face of changing product demands, $S = \{S_k\}_{k \in \mathcal{K}}$. However, solving \eqref{eqn:blending_prob} is highly challenging due to the presence of bilinearities in constraints \eqref{eqn:node_quality} and \eqref{eqn:quality_satis}. In fact, efficient solution of the blending problem has been the subject of much research due to its importance in the operation of petroleum refineries and waste-water treatment plants, among others  \citep{tawarmalani2002pooling, misener2009advances, gupte2015pooling}. To address this issue, we propose a solution that involves training a convex DFSOM of \eqref{eqn:blending_prob} offline for online use. In what follows, we show that our approach can significantly speed up the online solution of the blending problem while preserving solution accuracy, as measured against the \textit{global solutions} of \eqref{eqn:blending_prob}.

\subsubsection{Design and training of the DFSOM}
Our approach involves linearizing the nonconvex terms in \eqref{eqn:blending_prob}. For each $j \in \mathcal{J}$, $k \in \mathcal{K}$, and $w \in \mathcal{W}$, there is a bilinear term $q_{jw}x_{jk}$ in \eqref{eqn:blending_prob}, which we replace with the following approximation:
\begin{equation}
\sum\limits_{k' \in \mathcal{K}} \sum\limits_{\ell = 0}^{2} \left( p'_{jkwk'\ell} S_{k'}^{\ell} \, q_{jw} + q'_{jkwk'\ell} S_{k'}^{\ell} \, x_{jk} \right)
\end{equation}
where $p'_{jkwk'\ell}$ and $q'_{jkwk'\ell}$ are scalar parameters to be determined by solving \eqref{eqn:DFSM}. Similar to the previous two problems, we penalize the sum of absolute values of $p'_{jkwk'\ell}$ and $q'_{jkwk'\ell}$ in the objective function to induce sparsity. In this problem, we keep the objective function in the decision-focused surrogate problem the same as in the original problem.

We test the proposed decision-focused surrogate modeling approach on the blending network presented in Example 2 of Foulds et al.\cite{foulds1992bilinear} This network has four pooling nodes blending a total of six inlet streams that result in four final products; for the exact parameter values used in this case study, we refer to Foulds et al.\cite{foulds1992bilinear} To generate training data for \eqref{eqn:DFSM}, we solve \eqref{eqn:blending_prob} with $S_k$ values sampled from the uniform distribution $\mathcal{U}(100, 200)$ for all $k$ in $\mathcal{K}.$ A single data point here consists of the input vector $S$ and the corresponding optimal solution vector $(f, x, q).$ 

\subsubsection{Data efficiency of DFSOM in comparison to black-box optimization proxies}
\begin{figure}[h!]
    \centering
    \subfloat[Effect of training dataset size on the accuracy of DFSOMs.]{
        \label{fig:bp_data_eff}
        \begin{tikzpicture}[scale=0.9]
            \pgfplotstableread[col sep=comma]{Figures/pooling_pred_error.csv}{\csvdata}
            \begin{axis}[
                boxplot/draw direction = y,
                x axis line style = {opacity=100, thick},
                y axis line style = {thick},
                enlarge y limits,
                ymajorgrids,
                xtick = {1, 2, 3, 4, 5, 6},
                xticklabel style = {align=center},
                xticklabels = {$10$, $20$, $50$, $75$, $100$, $200$},
                xtick style = {}, 
                ylabel = {mean $\lvert z^* - z^{\mathrm{DFSOM}} \rvert $},
                ymax = 50,
                xlabel = {$|\mathcal{I}|$}
            ]
                \foreach \n in {1,...,6} {
                    \addplot+[boxplot = {every median/.style={red,ultra thick}}, draw=blue] table[y index=\n] {\csvdata};
                }
            \end{axis}
        \end{tikzpicture}} \hfil
    \subfloat[Effect of training dataset size on the accuracy of NNs.]{
        \label{fig:bp_nn_comparison}
        \begin{tikzpicture}[scale=0.7]
            \pgfplotstableread[col sep=comma]{Figures/E2E_440.csv}{\csvdatasmall}
            \pgfplotstableread[col sep=comma]{Figures/E2E_2622.csv}{\csvdatamed}
            \pgfplotstableread[col sep=comma]{Figures/E2E_4572.csv}{\csvdatalarge}
            \pgfplotstabletranspose\datatransposedsmall{\csvdatasmall} 
            \pgfplotstabletranspose\datatransposedmed{\csvdatamed} 
            \pgfplotstabletranspose\datatransposedlarge{\csvdatalarge} 
            \begin{axis}[
                boxplot/draw direction = y,
                x axis line style = {opacity=100, thick},
                y axis line style = {thick},
                width=\textwidth,
                height=7cm,
                enlarge y limits,
                ymajorgrids,
                xtick = {3, 4, 5, 6, 7, 8, 9, 10, 11, 12, 13},
                xticklabel style = {align=center},
                xticklabels = {$10$, $20$, $50$, $75$, $100$, $200$, $400$, $800$, ${1,600}$, ${3,200}$, ${4,000}$},
                xtick style = {}, 
                ylabel = {mean $\lvert z^* - z^{\mathrm{NN}} \rvert $},
                legend entries={nn-440,nn-2622,nn-4572},
                ymax = 225,
                xmin = 3,
                xmax = 13,
                enlarge x limits,
                xlabel = {$|\mathcal{I}|$}
            ]
                \addlegendimage{no markers, red}
                \addlegendimage{no markers, blue}
                \addlegendimage{no markers, black}
                \foreach \n in {3,...,13} {
                    \addplot+[boxplot = {every median/.style={red,ultra thick}, draw position=\n - 0.3, box extend = 0.2}, draw=red] table[y index=\n] {\datatransposedsmall};
                }
                \foreach \n in {3,...,13} {
                    \addplot+[boxplot = {every median/.style={blue,ultra thick}, draw position=\n, box extend = 0.2}, draw=blue] table[y index=\n] {\datatransposedmed};
                }
                \foreach \n in {3,...,13} {
                    \addplot+[boxplot = {every median/.style={black,ultra thick}, draw position=\n + 0.3, box extend = 0.2}, draw=black] table[y index=\n] {\datatransposedlarge};
                }
            \end{axis}
        \end{tikzpicture}}
    \caption{\label{fig:bp_surr_data_eff} The true optimal objective function value of \eqref{eqn:blending_prob} is denoted by $z^*$, whereas $z^{\mathrm{DFSOM}}$ and $z^{\mathrm{NN}}$ denote the objective values obtained by using the solutions generated by DFSOM and NN proxies, respectively. The box plots show the interquartile ranges of prediction error for the surrogate models obtained with ten different instances of randomly generated training data. In both figures, we use a separate test dataset of 100 samples to compute the prediction error.}
\end{figure}
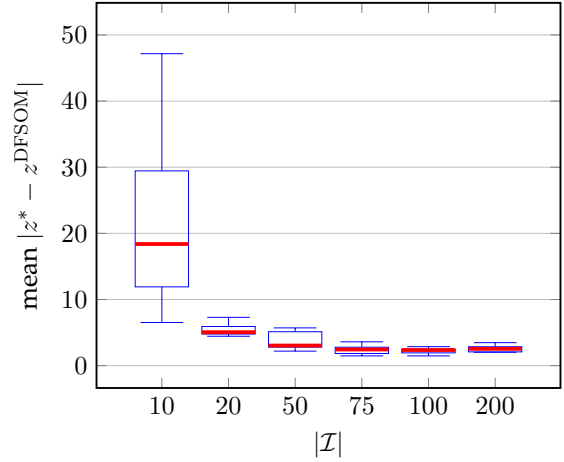
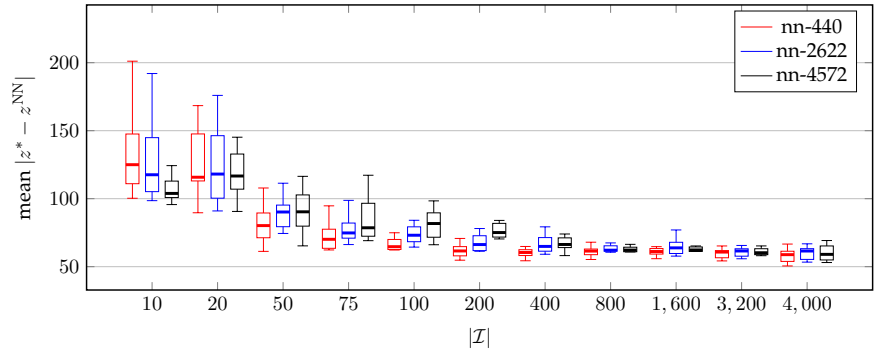

We examine the minimum amount of training data needed to create a high-quality surrogate. To do this, we solve \eqref{eqn:DFSM} using different sizes of $\mathcal{I}$, and then assess the prediction error on a separate test dataset of 100 points. Figure \ref{fig:bp_data_eff} displays the results of this analysis. Each box plot in the figure represents the interquartile range of prediction error for ten decision-focused surrogate models, each constructed using a different random training dataset of size $\lvert \mathcal{I} \rvert$. Remarkably, we find that, even for this high-dimensional problem, a low prediction error can be achieved with just 20 samples. Additionally, the prediction error achieves its minimum value with as few as 75 samples. These findings demonstrate the effectiveness of the proposed approach in developing surrogates with a small dataset, which is crucial since solving \eqref{eqn:blending_prob} to construct a large training dataset can be computationally burdensome.

The results in Figure \ref{fig:bp_data_eff} are in direct contrast to the ones presented in Figure \ref{fig:bp_nn_comparison} where we assess the accuracy of NN models that are trained to be optimization proxies for \eqref{eqn:blending_prob}. We train three different types of networks to produce optimal $f$ and $x$ values given inputs $S$. The main difference between these networks is their size, which allows them to have different numbers of learnable parameters (our DFSOMs have 288 learnable parameters); the key characteristics of these networks are summarized in Table \ref{tab:nn_sizes}. We train these NNs with different numbers of training data points ranging from 10 to 4,000. For each NN and $\lvert \mathcal{I} \rvert$ combination, we train ten different models, each with a different random training dataset. Similar to the DFSOM case, we evaluate the prediction accuracy of the NNs on a separate test dataset of 100 samples.

As shown in Figure \ref{fig:bp_nn_comparison}, even with 4,000 data points, the NNs exhibit significantly worse performance compared to the DFSOMs. We believe this is a direct consequence of the black-box nature of NNs, which does not allow for the incorporation of known correlations about the data, in contrast to the DFSOM. It is also important to note that the DFSOM output will always satisfy the crucial mass balances in \eqref{eqn:node_mb} and demand satisfaction constraints in \eqref{eqn:demand_satis}, whereas there is no trivial way of ensuring the same for the NN outputs. These findings further highlight the fact that the proposed decision-focused surrogate modeling approach can be a superior alternative to NNs when it comes to generating high-quality surrogates for computationally challenging optimization problems.

\begin{table}[ht]
\caption{Key characteristics of NNs trained as optimization proxies for \eqref{eqn:blending_prob}. We use ReLU activation functions for the hidden layers.}
\label{tab:nn_sizes}
\centering
\begin{tabular}{@{}cccc@{}}
\toprule
\multirow{2}{*}{NN name} & {number of} & {number of} & {number of} \\ 
& hidden layers & nodes/layer & learnable parameters \\ \midrule
nn-440 & 2 & 11 & 440 \\
nn-2622 & 4 & 25 & 2622 \\
nn-4572 & 7 & 25 & 4572 \\
\bottomrule
\end{tabular}
\end{table}

\subsubsection{Computational performance of DFSOMs for RTO}
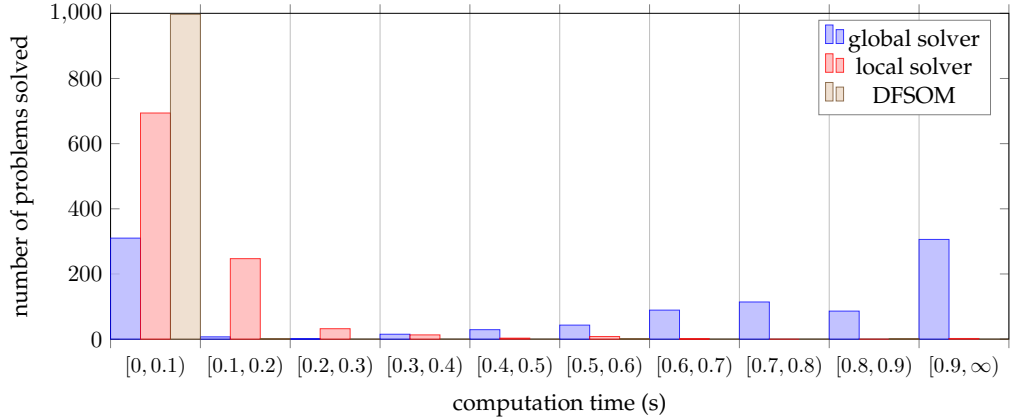
\begin{figure}[h!]
    \centering
\begin{tikzpicture}[scale = 0.8]
    \begin{axis}[
        ybar interval,
        ymax = 1000,
        ymin = 0,
        xmax = 1,
        xmin = 0,
        width=\textwidth,
        height=7cm,
        legend entries = {global solver, local solver, DFSOM},
        legend style = {draw=gray},
        ylabel = number of problems solved,
        xlabel = computation time (s),
        xtick={0, 0.1, 0.2, 0.3, 0.4, 0.5, 0.6, 0.7, 0.8, 0.9, 1},
        xticklabel style = {font=\small},
        xticklabels = {{$[0, 0.1)$}, {$[0.1, 0.2)$}, {$[0.2, 0.3)$}, {$[0.3, 0.4)$}, {$[0.4, 0.5)$}, {$[0.5, 0.6)$}, {$[0.6, 0.7)$}, {$[0.7, 0.8)$}, {$[0.8, 0.9)$}, {$[0.9, \infty)$}, 1}
    ]
        \addplot+ [opacity=0.8, hist={data=x, bins = 10, data min=0, data max = 1}] table [col sep = comma, x = global_true_time] {Figures/computation_time_stats.csv};
        \addplot+ [opacity=0.8, hist={data=x, bins = 10, data min=0, data max = 1}] table [col sep = comma, x = local_true_time] {Figures/computation_time_stats.csv};
        \addplot+ [opacity=0.8, hist={data=x, bins = 10, data min=0, data max = 1}] table [col sep = comma, x = surrogate_time] {Figures/computation_time_stats.csv};
    \end{axis}
\end{tikzpicture}
\caption{\label{fig:pooling_comp_stats} Computational performance of DFSOM in comparison to when \eqref{eqn:blending_prob} is solved with local and global optimization solvers.}
\end{figure}

We analyze the computational performance of the DFSOM in contrast to \eqref{eqn:blending_prob}. We solve 1,000 instances of a DFSOM and \eqref{eqn:blending_prob}, with each instance corresponding to different $S$ values following the same distribution as the training data. The optimization solver Gurobi is used to solve \eqref{eqn:DFSM}, while we solve each instance of \eqref{eqn:blending_prob} twice, once with the global solver Gurobi and again with the local solver IPOPT. A histogram in Figure \ref{fig:pooling_comp_stats} shows the distribution of computation times for the three cases. We find that, depending on $S$, the global solution of \eqref{eqn:blending_prob} can take more than 100 seconds to find, possibly making it intractable for online application. While the local solver provides a fast solution in most cases, the generated solutions generally suffer from large optimality gaps as shown in Table \ref{tab:blending_data}. IPOPT-generated solutions had a mean optimality gap of 61\%, which can be as high as 230\% in the worst case.

In contrast, DFSOM is highly effective in alleviating computational burden associated with the global solution of this challenging optimization problem. We find that all 1,000 instances of DFSOM solve in less than 0.1 seconds, and even the worst-case optimality gap is only 0.9\%. Therefore, in addition to being a superior surrogate modeling strategy, using decision-focused surrogate modeling can also be a better alternative to suboptimal local solvers.

\begin{table}[ht]
\caption{Optimality gap for solutions produced by DFSOM compared with when \eqref{eqn:blending_prob} is solved by IPOPT.}
\label{tab:blending_data}
\centering
\begin{tabular}{@{}cccc@{}}
\toprule
 & mean & maximum & minimum \\ \midrule
\multirow{2}{*}{$\frac{\lvert z^{\mathrm{local}} - z^{*} \rvert}{\lvert z^{*} \rvert}$}& \multirow{2}{*}{0.61} & \multirow{2}{*}{2.30} & \multirow{2}{*}{0} \\
& & & \\
\multirow{2}{*}{$\frac{\lvert z^{\mathrm{DFSOM}} - z^{*} \rvert} {\lvert z^{*} \rvert}$}& \multirow{2}{*}{0.002} & \multirow{2}{*}{0.009} & \multirow{2}{*}{0} \\
& & & \\
\bottomrule
\end{tabular}
\end{table}

\begin{remark}
In all three numerical case studies, we employ polynomial functions of inputs to represent the nonlinear coefficient models. For instance, in \eqref{eqn:ratelaw_surr}, $a(A_i)$, $b(A_i)$, and $c(A_i)$ are modeled as cubic polynomials in $A_i$. The choice of polynomials stems from their linearity in learnable parameters, thus preserving the multiconvex nature of \eqref{eqn:KKT-EPR}. While the proposed framework can accommodate more complex functions like neural networks with nonlinear activation functions, this would preclude the use of BCD for solving the ensuing learning problem, potentially complicating the solution of larger instances. Nonetheless, our computational findings, particularly depicted in Figure \ref{fig:cstr_feasible_region}, demonstrate that the proposed framework can achieve highly accurate decision predictions from simpler models. This can be attributed to the framework's focus on learning only the optimal solution space of an optimization problem, rather than attempting to model the entire feasible solution space.
\end{remark}

\section{Conclusions} \label{sec:conclude}
In this work, we introduced a novel decision-focused surrogate modeling approach for the online solution of computationally challenging nonlinear optimization problems. Our data-driven framework produces a surrogate optimization problem that minimizes the decision prediction error on a training dataset containing optimal solutions of the original optimization problem corresponding to different model inputs. We showed that the learning problem can be formulated and solved as a data-driven inverse optimization problem. Through three computational case studies, we demonstrated that:

\begin{enumerate}
\item Our approach produces high-quality surrogates with much simpler surrogate representations of the feasible regions of the original problem compared to traditional methods that involve optimization with embedded surrogate models. This key benefit arises from the decision-focused nature of our approach as it does not seek to learn the entire feasible solution space.
\item Simpler models of the feasible region lead to convex surrogate optimization problems, which obviates the need for expensive global solvers while still generating solutions that are close to globally optimal solutions.
\item Compared to black-box models used as optimization proxies, our approach is significantly more data-efficient, allowing the user to retain a large part of the original optimization problem that does not contribute to the problem's nonconvexity. 
\item Our framework produces an optimization problem as the resulting surrogate model, making it easier to incorporate essential system constraints as hard constraints, which is typically not straightforward with black-box optimization proxies based on, for example, neural networks.
\end{enumerate}

In summary, our decision-focused surrogate modeling paradigm presents a promising new avenue for the solution of time-critical optimization problems, offering both higher quality and more efficient surrogate optimization problems in comparison to traditional surrogate modeling methods. In a follow-up work \cite{gupta2022decision}, we extended the proposed framework to include a mechanism that minimizes potential infeasibility (with respect to the original problem) of the DFSOM's solutions. Several important directions for future work still remain, including the nonparametric construction of convex surrogate models to replace the nonconvex parts of the original problem, and the development of adaptive sampling algorithms to further improve data efficiency.

\section*{Data Availability and Reproducibility Statement}
The computer code and datasets that support the findings of this study are available under the MIT license at https://github.com/ddolab/DecFocSurrMod. Specifically, plots in Figure 2 can be regenerated by running the Julia code in the ''cstr" folder. Steps to reproduce Figures 3 and 4 are shown in the \texttt{plots.ipynb} Jupyter notebook in the same folder. The Julia and Python code to reproduce results in Figure 6 is available in the ''heat\_exchanger" folder. Finally, the input files to regenerate the results in Figures 7 and 8 as well as Table 2 are available in the ''blending\_network" folder. The computer code in the GitHub repository also makes available the exact algorithmic hyperparameter settings used for the results shown in this paper. 

\section*{Acknowledgements}
The authors thank Sayandeep Biswas and Joshua Larson, who were at the time undergraduate students at the University of Minnesota, for conducting preliminary computational analysis for this work. The authors gratefully acknowledge the financial support from the National Science Foundation under Grant \#2044077 as well as the Minnesota Supercomputing Institute (MSI) at the University of Minnesota for providing resources that contributed to the research results reported in this paper. R.G. acknowledges financial support from a departmental fellowship sponsored by 3M and a Doctoral Dissertation Fellowship from the University of Minnesota.

\bibliographystyle{apalike}
\bibliography{library,responses}

\end{document}